\documentclass[12pt, leqno]{amsart}
\usepackage[all]{xy}
\usepackage{amsfonts,amsmath,amssymb,amscd}
\usepackage{tabularx,ragged2e,booktabs,caption}

\newtheorem{theorem}{Theorem}

\newtheorem{lemma}[subsection]{{\bf Lemma}}

\newtheorem{remark}[subsection]{Remark}

\newcommand{\al}{\alpha}

\newcommand{\Z}{\mbox{$\mathbb Z$}}

\begin{document}
	\title{ Extension of Irreducibility results on Generalised Laguerre Polynomials $L_n^{(-1-n-s)}(x)$} 
\author[Nair]{Saranya G. Nair}
\address{Department of Mathematics\\
BITS Pilani, K K Birla Goa Campus, Goa- 403726}
\email{saranyan@goa.bits-pilani.ac.in}
\author[Shorey]{T. N. Shorey}
\address{National Institute of Advanced Studies, IISc Campus\\
Bangalore, 560012}
\email{shorey@math.iitb.ac.in}

\begin{abstract} We consider the irreducibility of Generalised Laguerre Polynomials
 for negative integral values given by $L_n^{(-1-n-s)}(x)=\displaystyle\sum_{j=0}^{n}\binom{n-j+s}{n-j}\frac{x^j}{j!}.$ For different values of $s,$ this family gives polynomials which are of great interest. It was proved earlier that for $0 \leq s \leq 60,$ these polynomials are irreducible over $\mathbb{Q}.$ In this paper we improve this result upto $s \leq 88.$
\end{abstract}

\maketitle
\pagenumbering{arabic}
\pagestyle{myheadings}
\markright{Extension of Irreducibility results on Generalised Laguerre Polynomials $L_n^{(-1-n-s)}(x)$}
\markleft{ Nair and Shorey}
\section{{\bf Introduction}}

For a positive integer $n$ and real number $\alpha,$ the Generalised Laguerre Polynomial (GLP) is defined as
\begin{align}\label{lag}
L_n^{(\alpha)}(x)=\displaystyle\sum_{j=0}^{n}\frac{(n+\al)(n-1+\al)\cdots (j+1+\al)}{j!(n-j)!}(-x)^j.
\end{align}
These polynomials were discovered around 1880 and they have been extensively studied in various branches of mathematics and mathematical physics. The algebraic properties of GLP were first studied by Schur \cite{Sch1},\cite{Sch2} where he established the irreducibility of $L_n^{(\alpha)}(x)$ for $\alpha \in \{0,1,-n-1\},$ gave a formula for the discriminant $\Delta_n^{(\alpha)}$ of $ \mathcal{L}_n^{(\alpha)}(x)=n!L_n^{(\alpha)}(x)$ by \begin{align*}
	\Delta_n^{(\alpha)}=\displaystyle\prod_{j=1}^{n}j^j(\alpha+j)^{j-1}
\end{align*}
and calculated their assosciated Galois groups. For an account of results obtained on GLP, we refer to \cite{Haj},\cite{NaSh}.

Let $f(x) \in \mathbb{Q}[x]$ with deg $f=n$. By irreducibility of a polynomial, we shall always mean its irreducibility over $\mathbb{Q}.$ We observe that if a polynomial of degree $n$ has a factor of degree $k <n,$ then it has a factor of degree $n-k.$ {\it Therefore given a polynomial of degree $n$, we always consider factors of degree $k$ where $1 \leq k \leq \frac{n}{2}$}. If the argument $\alpha$ of \eqref{lag} is a negative integer, we see that the constant term of $L_n^{(\alpha)}(x)$ vanishes if and only if $n \geq |\alpha|=-\alpha$ and then $L_n^{(\alpha)}(x)$ is reducible. Therefore we assume that $\alpha \leq -n-1.$ We write $\alpha=-n-s-1$ where $s$ is a non-negative integer. We have 
\begin{align}
	L_n^{(-n-s-1)}(x)=\displaystyle\sum_{j=0}^{n}(-1)^n\frac{(n+s-j)!}{(n-j)! s!}\frac{x^j}{j!}.
\end{align}
Borrowing the notation from \cite{SinSho}, we consider the following polynomial 
\begin{align*}
	g(x):=g(x,n,s)=(-1)^nL_n^{(-n-s-1)}=\displaystyle\sum_{j=0}^{n}\binom{n+s-j}{n-j}\frac{x^j}{j!}=\displaystyle\sum_{j=0}^{n}b_j\frac{x^j}{j!}
\end{align*}
where $b_j=\binom{n+s-j}{n-j}$ for $0 \leq j \leq n.$ Thus $b_n=1,b_0=\binom{n+s}{s}=\frac{(n+1)\cdots (n+s)}{s!}$. We observe that $g(x)$ is irreducible if and only if $L_n^{(-n-s-1)}(x)$ is irreducible. The aim of this paper is to discuss the irreducibility of $g(x)$. We consider more general polynomial
\begin{align*}
	G(x):=G(x,n,s)=\displaystyle\sum_{j=0}^{n}a_jb_j\frac{x^j}{j!}
\end{align*}
such that $a_j \in \mathbb{Z}$ for $0 \leq j \leq n$ with $|a_0|=|a_n|=1.$ If $a_j=1$ for $0 \leq j \leq n,$ we have $G(x)=g(x).$ We write
\begin{align*}
	g_1(x):=n! g(x) = n! \displaystyle\sum_{j=0}^{n}\binom{n+s-j}{n-j}\frac{x^j}{j!}
\end{align*}  and 
\begin{align*}
	 G_1(x)=n! G(x) 
\end{align*}so that $g_1$ and $G_1$ are monic polynomials with integer coefficients of degree $n.$ The irreducibility of $g_1(x)$ and $G_1(x)$ implies the irreducibility of $g(x)$ and $G(x)$ respectively.  We begin with the following result by Sinha and Shorey \cite{SinSho} on $G(x)$. 

\begin{lemma}\label{Lemma 1}
Let $s \leq 92.$ Then $G(x)=G(x,n,s)$ has no factor of degree $k \geq 2$ except when $(n,k,s) \in \{ (4,2,7),(4,2,23),(9,2,19),(9,2,47),(16,2,14),(16,2,34),(16,2,89),\\(9,3,47), (16,3,19),(10,5,4)\}. $
\end{lemma}
 We re-state Lemma 1.1 as follows:\\
  Let $s \leq 92.$ Assume that $G(x)$ has a factor of degree greater than or equal to $2$. Then
$(n,k,s) \in \{ (4,2,7),(4,2,23),(9,2,19),(9,2,47),(16,2,14),(16,2,34),(16,2,89),(9,3,47),\\ (16,3,19),(10,5,4)\}.$ We check that $g(x)$ is irreducible when $n \in \{ 4,9,10,16\}.$ Therefore by Lemma \ref{Lemma 1} with $G(x)=g(x)$ we have 
\begin{lemma}\label{Lemma2}
	Let $n \geq 3$ and $s \leq 92$. Then $g_1(x)$ is either irreducible or linear factor times an irreducible polynomial. 
\end{lemma}
The irreducibility of $g_1(x)$ was proved by Schur \cite{Sch1} for $s=0,$ by Hajir \cite{Haj06} for $s=1$, by Sell \cite{Sell} for $s=2$ and by Hajir \cite{Haj} for $3 \leq s \leq 8.$ We shall prove 
\begin{theorem}\label{thm1}
	$g_1(x)$ is irreducible for $9 \leq s \leq 88.$ 
\end{theorem}  
Nair and Shorey \cite{NaSh15b} and Jindal, Laishram and Sarma \cite{jin} already proved the irreducibility of $g_1(x)$ in the range of $9 \leq s \leq 22$ and $23 \leq s \leq 60$ respectively. But our proof of Theorem 1 is new. The proofs of \cite{NaSh15b} and \cite{jin} depend on the method of Hajir in \cite{Haj} whereas our proof depends on Lemma \ref{Lemma2} which is a direct consequence of Lemma \ref{Lemma 1}. We could not cover the cases of $88 \leq s \leq 92$ due to computational limitations.

 At this point, we pause with a digression concerning notation.
In order to dispell possible confusion in the reader concerning notation, it is worth pointing out that the notation $L_n^{<s>}(x)$ (resp., $\mathcal{L}_n^{<s>}(x)$) was used by the authors in \cite{Haj}, \cite{jin} and \cite{NaSh} to denote the polynomial $g(x, n, s)$ (resp, $g_1(x, n, s)$).

\section{Preliminaries}
 From now onwards we shall assume that $s \geq 9.$  For a real number $\alpha,$ we write $\left[ \alpha\right]$ to be the largest integer not exceeding $\alpha.$
 
  Let $f(x) =\displaystyle\sum_{j=0}^{m}d_jx^j \in \Z[x]$ with $d_0d_m \neq 0$ and let $p$ be a prime. For an integer $x,$ let $\nu(x)=\nu_p(x)$ be the highest power of $p$ dividing $x$ and we write $\nu(0)=\infty.$ Let S be the following set of
points in the extended plane
$$S =\{(0,\nu(d_m)),(1,\nu(d_{m-1})),(2,\nu(d_{m-2})),\ldots,(m,\nu(d_0))\}.$$
Consider the lower edges along the convex hull of these points. The left most endpoint
is $(0,\nu(d_m))$ and the right most endpoint is $(m,\nu(d_0))$. The endpoints of each edge
belong to S and the slopes of the edges increase strictly from left to right. The polygonal path formed by these edges is called the Newton polygon
of $ f(x)$ with respect to the prime $p$ and we denote it by $NP_p(f)$. The endpoints of the
edges of $NP_p(f)$ are called the vertices of $ NP_p(f)$.
We begin with a very useful result, due to Filaseta \cite{Fil}, giving a criterion on the factorisation of a polynomial in terms of the maximum slope of the edges of its Newton polygon. 
\begin{lemma}\label{newton1}
	Let $l, k,m$ be integers with $m \geq 2k > 2l \geq 0$. Suppose $h(x) =\displaystyle\sum_{j=0}^{m}b_jx^j \in
	\Z[x] $ and $p$ be a prime such that $ p \nmid b_m$ and $p\mid b_j $ for $0 \leq j \leq m-l-1$ and the
	right most edge of $ NP_p(h)$ has slope $ <\frac{1}{k}$. Then for any integers $a_0, a_1,\ldots, a_m$
	with $p\nmid a_0a_m$, the polynomial $f(x) =\displaystyle\sum_{j=0}^{m} a_jb_jx^j$ cannot have a factor with degree in $[l + 1, k]$.
\end{lemma}
The next result is Lemma 4.2 of \cite{ShTi10} with $a=0.$
\begin{lemma}\label{ShTi}
	Let $a_0,a_1,\cdots,a_n$ denote arbitrary integers and $$h(x)=\displaystyle\sum_{j=0}^{n}a_j\frac{x^j}{j!}.$$ Assume that $h(x)$ has a factor of degree $k \geq 1.$ Suppose that there
	exists a prime $p > k $ such that $p$ divides $n(n-1) \cdots (n-k+1).$ Then $p$ divides $a_0a_n.$
\end{lemma}
For a positive integer $l$ and a prime $p,$ let $\nu_p(l)$ be the maximal power of $p$ dividing $l.$ 

\begin{lemma}\label{order1}
	Let $p$ be a prime. For any integer $l\geq 1$, write $l$ in base $p$ as $l=l_tp^t+l_{t-1}p^{t-1}+\dots+l_1p+l_0$ where $0\leq l_i\leq p-1$ for $0\leq i \leq t$ and $l_t >0$. Then
	\begin{align*}
		\nu_p(l!)=\frac{l-\sigma_p(l)}{p-1}
	\end{align*}
	where $\sigma_p(l)=l_t+l_{t-1}+\dots+l_1+l_0$.
\end{lemma}
This is due to Legendre. For a proof, see \cite[Ch.17, p 263]{Hasse}. As a consequence we have
\begin{align}\label{eq2}
	\nu_p\left(\binom{m}{t}\right)= \frac{\sigma_p(t)+\sigma_p(m-t)-\sigma_p(m)}{p-1}.
\end{align}
	\begin{lemma}\label{Lemma 6}
		Assume that $g_1(x)$ is a linear factor times an irreducible polynomial. Let $p$ be a prime dividing $n$ and $ s < p^2.$ Then \begin{align*}
			d+ \left[\frac{s}{p}\right] \geq p
		\end{align*}
	where $d \equiv \frac{n}{p} (\mod p)$ for $ 1 \leq d <p.$ 
	\end{lemma}
The assertion of Lemma \ref{Lemma 6} was proved in \cite[Corollary 3.2]{jin} under the assumption of $p$ dividing $n_1$ where 

\begin{align}\label{defn of n}
	n=n_0\cdot n_1 \ {\rm with} \ \gcd(n_0,n_1)=1
\end{align}
and
\begin{align}\label{n_1}
	n_1=\displaystyle \prod_{p| \gcd (n,\binom{n+s}{s})}p^{\text{ord}_p(n)}.
\end{align} Therefore $n_0$ is the largest divisor of $n$ which is coprime to $\binom{n+s}{s}.$ Thus the assumption $p$ dividing $n_1$ in \cite[Corollary 3.2]{jin} is replaced by $p$ dividing $n$ in Lemma \ref{Lemma 6} when $g_1(x)$ is linear factor times an irreducible polynomial.

\begin{proof}
	We apply Lemma \ref{ShTi} with $h(x)=g_1(x)$ and $k=1$ to conclude that \begin{align} \label{eq1}
		p|\frac{(n+1)\cdots (n+s)}{s!}.
	\end{align}
If $\nu_p(n)>s,$ then $\nu_p\left( \frac{(n+1)\cdots (n+s)}{s!}\right) \leq \nu_p(\frac{s!}{s!})  \leq 0$ contradicting \eqref{eq1}. Therefore $\nu_p(n) \leq s < p^2$ and hence $ 1 \leq \nu_p(n) <2.$ We write,
\begin{align*}
	n=pD \ \text{where} \ \gcd(D,p)=1
\end{align*} and $s=ps_1+s_0$ where $1 \leq s_1 <p, 0 \leq s_0 <p$. Then $n+s=p(D+s_1)+s_0$ which implies that $\sigma_p(n+s)=\sigma_p(D+s_1)+s_0$. Now we argue as in \cite[Lemma 3.1]{jin} for deriving from \eqref{eq2} and \eqref{eq1} that 
\begin{align*}
1 \leq \nu_p\left(\binom{n+s}{s}\right)= &\frac{\sigma_p(n)+\sigma_p(s)-\sigma_p(n+s)}{p-1}\\
=&	\frac{\sigma_p(D)+s_1+s_0+-\sigma_p(D+s_1)-s_0}{p-1}\\
=&\frac{\sigma_p(D)+s_1-\sigma_p(D+s_1)}{p-1}\\
=& \nu_p\left(\binom{D+s_1}{s_1}\right)\\
=& \nu_p \left(\frac{(D+1)\cdots (D+s_1)}{s_1!}\right)\\
=& \nu_p((D+1) \cdots (D+s_1))\\
=& \nu_p(D+j) \ \text{for \ exactly one} \ j \ \text{with} \ 1 \leq j \leq s_1 
\end{align*}
since $s_1 <p.$ Hence $D+\left[ \frac{s}{p}\right]=D+s_1 \geq D+j \equiv 0 $ (mod $p$). Since $\gcd(D,p)=1$, we observe from \eqref{eq2} that $D=\frac{n}{p} \equiv d$ (mod $p$) where $ 1 \leq d <p$ by looking at $p-$ adic representation of $n.$ Hence $ d+ \left[ \frac{s}{p}\right] \geq p$ where $d \equiv \frac{n}{p}$ (mod $p$) for $1 \leq d <p.$ 
\end{proof}
Next, we prove
\begin{lemma}\label{Lemma 7}
	Assume that $g_1(x)$ is linear factor times an irreducible polynomial. Then for $n \leq 127$ and $s \leq 103,$ $g_1(x)$ is irreducible.
	\end{lemma}
	This result is proved in \cite[Lemma 2.10]{NaSh15b} without the assumption that $g_1(x)$ is linear factor times an irreducible polynomial. But Lemma \ref{Lemma 7} suffices for our purpose. We give here a proof of this particular case, as it involves considerably less computations.
	\begin{proof}
		Let $n \leq 127$ and $s \leq 103.$ Since $g_1(x)$ is linear factor times an irreducible polynomial, we see that $n_0=1.$ Assume that $n \geq s+2.$ Then $\max( \frac{n+s}{2},n-1)=n-1$ and we derive from \cite[Lemma 4.1]{Haj} that $g_1(x)$ is irreducible if $n$ is prime. Now we check that $g_1(x)$ is irreducible for all pairs $(n,s)$ with $n$ composite and $n \geq s+2.$ Next we assume that $n <s+2.$ Then $\max( \frac{n+s}{2},n-1)=\frac{n+s}{2}$. We determine all pairs $(n,s)$ such that $n <s+2$ and there exists a prime $p$ satisfying $\frac{n+s}{2} < p \leq n$. Then $g_1(x)$ is irreducible for all these pairs $(n,s)$ by \cite[Lemma 4.1]{Haj}. For the remaining pairs $(n,s)$ with $n <s+2,$ we check that $g_1(x)$ is irreducible. 
		\end{proof}
	We close this section by stating the following result which is Lemma 3.6 from \cite{jin} with $i=1$.
	\begin{lemma} \label{Lemma factor1}
		Let $p|n(s+1)$ and $\nu_p\left(\binom{n+s}{s}\right)=u.$ Then $g_1(x)$ cannot have a factor of degree $1$ if any of the following condition holds: 
	\begin{itemize}
		\item[(i)] $u=0$
		\item[(ii)] $u >0$, $ p>2$ and $
		\max \left\{\frac{u+1}{p}, \frac{\nu_p(n+s-z_0)-\nu_p(n)}{z_0+1}\right\} < 1$,
 where $z_0 \equiv n+s$ (mod $p$) with $1 \leq z_0 <p.$
 	\end{itemize}
		\end{lemma}
	
	\section{Proof of Theorem \ref{thm1}}
	For $c> 1$ and $s \geq c^2$, we consider the following set given by
	\begin{align*}
		H_{s,c}=\{n \in \mathbb{N}, n >127 \ \text{and for } p |n, p^{\nu_p(n)} \leq s \ \text{and if} \ p > sc^{-1} \ \text{then} \ d+\left[\frac{s}{p}\right] \geq p \}
	\end{align*}
where $1 \leq d <p$ and $d \equiv \frac{n}{p}$ (mod $p$). Since $p \geq sc^{-1} \geq \sqrt{s},$ we derive from Lemmas \ref{Lemma 6} and \ref{Lemma 7} that it suffices to prove the irreducibility of $g_1(x)=g_1(x,n,s)$ with $n \in H_{s,c}.$ We partition $H_{s,c}$ as $H_{s,c,1}$ and $H_{s,c,2}$ given by 
\begin{align*}
	H_{s,c,1}=\{ n \in H_{s,c} \ \text{such that} \ P(n) \leq [sc^{-1}]\}
\end{align*}
and 
\begin{align*}
	H_{s,c,2}=\{ n \in H_{s,c} \ \text{such that} \ P(n) > [sc^{-1}]\}.
\end{align*}
Let $9 \leq s \leq 88$. By taking $ c \in \{3,3.42,5.5,7.7\}$ we compute $H_{s,c,1}$ and $H_{s,c,2}$ and hence $H_{s,c}$ for $ 9 \leq s \leq 88.$ We give some details regarding the computations of $H_{s,c}$ below. For example, for $s=80,$ and $ c=7.7$, the cardinality of $H_{s,c}$, $|H_{s,c}|=1538$ and for $s=85,$ and $ c=7.7$, the cardinality of $H_{s,c}$, $|H_{s,c}|=2466.$ The following table gives the $c$ values which we have chosen for each $s$ to compute the set $H_{s,c}$.
\begin{center}

\begin{tabular}{|c|c|}
	\hline
	$s$ & $c$ \\
	\hline
	$9 \leq s \leq 11$ & $3$ \\
	\hline
$12 \leq s \leq 35$ & $3.42$ \\
\hline
$36 \leq s \leq 60$ & $5.5$ \\
\hline	
 $61 \leq s \leq 88$ & $7.7$\\
 \hline
\end{tabular}
\end{center}
For each $n \in H_{s,c}$, we apply Lemma \ref{Lemma factor1} to derive that $g_1(x)=g(x,n,s)$ is irreducible except for $(n,s) \in T$ where 
\begin{align*}
	T=&\{(272,17),(144,21),(144,23),(144,25),(144,26),(312,26),(600,26),(216,29),(216,31),\\& (720,31),(240,35), (1440,35),(288,40),(288,41),(216,42),(216,44),(216,47),(288,47),\\& (288,48),(216,49),(144,51),(288,51),(144,53),(216,53),(288,53),(4320,55),(216,59),\\&(216, 63), (288, 63),(432, 63), (672, 63),(180,71),(192,71),(216,71),(216,79),(576,79),\\&(144, 80), (192, 80),(216, 80), (320, 80), (432, 80), (576, 
			80), (720, 80), (4320, 80)\}
\end{align*}
We use IrreduicibilityQ command in Mathematica to check that $g_1(x)$ is irreducible for all these values of $(n,s) \in T.$ 
\begin{remark}
	It is not necessary to compute $H_{s,c}$ for all values of $s$. For a fixed $c$, if $[sc^{-1}]=[(s+1)c^{-1}]$ and $P(s+1) \leq [sc^{-1}]$, then $H_{s,c}=H_{s+1,c}$. The assertion follows from the definitions of $H_{s,c,1}$ and $H_{s,c,2}$. Therefore $H_{s,c}=H_{s+1,c}$ for $s \in \{19,27,29,34\}$ with $c=3.42$, $s \in \{39,41,47,49,53,55,59\}$ with $c=5.5$ and $s \in \{62,69,71,74,79,83,87\}$ with $c=7.7$.
		\end{remark}
	\begin{remark}
		 We can take $c=5.5$ or $c=7.7$ according as $s$ at least $36$ or $s \geq 60$ respectively. But we cannot take $c$ more than $3.42\frac{s}{s-1}$ without sharpening \cite[Corollary 6]{NaSh15b}. Consequently we get sets $H_{s,c}$ of smaller size when $s \geq 36$ and this reduces the computations.
\end{remark}
\subsection*{Acknowledgements} The computations for calculating $H_{s,c}$  have been carried out by using SAGE,
and all other computations have been done in MATHEMATICA. We thank
Ankita Jindal 
for helping us with the computations. The second author is supported by INSA during the period of this work.
  
\end{document}